\documentclass[12pt,authoryear]{article}
\usepackage{amsfonts,amsthm,amsmath,mathrsfs,natbib} 
\usepackage{amssymb}
\usepackage{graphicx}
\usepackage{setspace} 
\usepackage{multirow}
\usepackage{xcolor}
\usepackage[utf8]{inputenc}
\usepackage[affil-it]{authblk}

\newtheorem{theorem}{Theorem} 

\theoremstyle{definition}
\newtheorem{remark}[theorem]{Remark}
\newtheorem{example}[theorem]{Example}

\def\BB{ \mathfrak{B} }

\newcommand{\R}{\mathbb{R}}
\newcommand{\N}{\mathbb{N}}
\newcommand{\sech}{{\rm sech\,}}
\newcommand{\norm}[1]{\lVert#1 \rVert}
\newcommand{\abs}[1]{\lvert#1 \rvert}

\allowdisplaybreaks

\begin{document}

\title{Mixture representations of noncentral distributions}

\author{Ludwig Baringhaus and Rudolf Gr{\"u}bel}

\affil{
  Leibniz Universit\"at Hannover\\
  Institut f\"ur Mathematische Stochastik\\
  Postfach 60 09, D-30060 Hannover, Germany\\
  lbaring@stochastik.uni-hannover.de\\
  rgrubel@stochastik.uni-hannover.de
}

\date{}

\maketitle

\begin{abstract} 
	With any symmetric distribution $\mu$ on the real line we may associate a parametric 
	family of noncentral distributions as the distributions of $(X+\delta)^2$,  $\delta\not=0$,
	where $X$ is a random variable with  distribution $\mu$. 
	The classical case arises if $\mu$ is the standard normal 
	distribution,  leading to the noncentral chi-squared distributions. It is well-known that 
	these may be written as Poisson  mixtures of the central chi-squared distributions with
	odd degrees of freedom. We obtain such mixture representations for the logistic distribution
	and for the hyperbolic secant distribution. We also derive  alternative representations for 
	chi-squared distributions and relate these to representations of the Poisson family.
	While such questions originated in parametric statistics they also appear 
	in the context of the generalized second Ray-Knight theorem, which connects Gaussian 
	processes and local times of Markov processes.       



\end{abstract}

\section{Introduction}\label{sec:intro}

Many of the classical parametric distributions,  such as the normal and Poisson families,
arise in the context of limit theorems, whereas others, such as the exponential distributions, 
are characterized by certain properties or by invariance under specific transformations. 
In view of this origin direct distributional relations between these families are 
often found to be surprising. A particularly intriguing example is the following, which is 
also the starting point for the present paper:  
\begin{equation}\label{eq:nzchiqua0}
(X + \delta)^2\; =_{\mathcal D} \;X^2+2\,\sum_{j=1}^N E_j.
\end{equation}
Here $\delta\neq 0$ is a real number, the random variables $N,X,E_1,E_2,\dots$ are independent,
$X$ has the standard normal distribution, $E_1,E_2,\dots$ are exponentially
distributed with mean 1, $N$ has the Poisson distribution with parameter $\delta^2/2$, and
`$=_{\mathcal D}$' denotes equality in distribution.

In statistics, the distribution in \eqref{eq:nzchiqua0} is known as 
the noncentral chi-squared distribution with one degree of freedom and noncentrality parameter 
$\delta^2$, and commonly denoted by $\chi^2_1(\delta^2)$. We recall that the
central chi-squared distributions $\chi_k^2$ with $k$ the degree of freedom, 
$k\in\N$,
the set of natural numbers not including 0, and the exponential distributions are both 
subfamilies of the family $\Gamma(\alpha,\lambda)$, $\alpha,\lambda>0$, of gamma distributions.
Specifically, $\chi_k^2=\Gamma(k/2,1/2)$ and $\Gamma(1,\lambda)$ is the 
exponential distribution with mean $1/\lambda$. In view of the convolution property of the 
gamma family,
\begin{equation*}
 \Gamma(\alpha,\lambda) \star \Gamma(\beta,\lambda) = \Gamma(\alpha+\beta,\lambda)
      \quad\text{for all } \alpha,\beta,\lambda >0,
\end{equation*}
we may thus rewrite \eqref{eq:nzchiqua0} as a mixture representation,
\begin{equation}\label{eq:chimix}
  \chi^2_1(\delta^2) \;=\; e^{-\delta^2/2} \sum_{n=0}^\infty \frac{\delta^{2n}}{2^n n!}\, \chi^2_{2n+1}
  	 \quad\text{for all } \delta\not=0,
\end{equation}
and indeed, it often appears in this form; see e.g.~pp.382,458 in~\cite{JKB}.
Interestingly, the equivalent relations \eqref{eq:nzchiqua0} and \eqref{eq:chimix}
appear under two quite different circumstances: In statistics in connection with the power of statistical tests,
see~\cite{liese}, and in probability theory in connection with the local times of Markov processes, 
see Lemma~6.33 in~\cite{moper}.

Note that the base family $\{\chi^2_k:\, k\in\N\}$ of mixture components is the same for all 
noncentrality parameters, so 
that~\eqref{eq:chimix} may be regarded as separating the two variables, the argument
of the probability measure and its parameter. A third way to write the representation is in terms 
of random variables $Y_k$ with distribution $\chi_{2k+1}^2$, $k\in\N_0$, 
the set of natural numbers including 0, as 
\begin{equation}\label{eq:repr3}
	(X+\delta)^2 =_{\mathcal D} Y_N,
\end{equation}
where the random variable $N$ is independent of the $Y$-variables and Poisson distributed with parameter $\delta^2/2$.

In the present paper we obtain similar results for two other families of noncentral distributions. 
For example, if $X$ has the hyperbolic secant distribution then, with $\phi:[0,\infty)\to [0,1)$,
$\phi(x)=2/(e^{2x} + e^{-2x})$,
\begin{equation}\label{eq:hypsec0}
    (X + \delta)^2\; =_{\mathcal D} \; \phi^{-1}\biggl(\frac{X_0^2}{X_0^2 +X_1^2 + 2\sum_{j=1}^{N} E_j}\biggr)^2.
\end{equation} 
Here, similar to~\eqref{eq:nzchiqua0}, the random variables $X_0,X_1,N,E_1,E_2,\ldots$   
are independent, $X_0$ and $X_1$ are standard normal, $E_1,E_2,\ldots$ are exponentially distributed with mean~1, 
but $N$ now has the negative binomial distribution with parameters $1/2$ and 
$2 /(e^{\delta}+e^{-\delta})^{2}$.

We also reconsider the classical case, where $\mu$ is the standard normal distribution. 
We find alternative representations for $\chi^2_1(\delta^2)$ and relate
these to mixture representations of the Poisson distributions. 

Section~\ref{sec:results} contains the main results, proofs are given in Section \ref{sec:proofs}. 
In Section~\ref{sec:comments} we collect various   
remarks; in particular, we expand on the connections to tests and local times mentioned above. 

\section{Main results}\label{sec:results}

In the first two subsections we obtain mixture representations for the noncentral distributions
associated with two non-normal distributions that are symmetric about 0. In the third 
subsection we return to the normal distribution and construct a general family of representations that
contains \eqref{eq:nzchiqua0} as a special case. 
These turn out to be related to mixture representations of the  Poisson distributions.

\subsection{The logistic distribution}\label{subsec:logistic}
The distribution function $F$ of the logistic distribution  is given by
\begin{equation}\label{eq:logDF}
	F(x) = \frac{1}{1+e^{-x}},\quad x\in\R.
\end{equation}
In view of $F(x)=1-F(-x)$ this distribution is symmetric about 0.
It will be convenient to rescale the noncentrality
parameter via
\begin{equation}\label{eq:delta_to_theta}
	\theta = \theta(\delta) := F(\delta) = \frac{e^\delta}{1+e^\delta}.
\end{equation}   
We consider two types of noncentrality. 

\begin{theorem}\label{thm:logistic}
	Suppose that $X$ has the distribution function $F$ given in \eqref{eq:logDF}. Then, for
	all $\delta\not=0$ and with $\theta$ as in \eqref{eq:delta_to_theta},
	\begin{align}\label{eq:reprlog1}
		|X+\delta|\ &=_{\mathcal D} \ \sum_{j=1}^N \frac{1}{j}\, E_j,\\
		\label{eq:reprlog2}
		(X+\delta)^2\ &=_{\mathcal D} \ \max\{ E_1^2, E_2^2,\ldots,E_N^2\},
	\end{align}
	where $N,E_1,E_2,\ldots$ are independent random variables, all $E_j$ are exponentially 
	distributed with mean 1, and
	\begin{equation}\label{eq:sumlog}
		P(N=n) \; = \; \theta\,(1-\theta)^n + (1-\theta)\,\theta^n\quad\text{for all }n\in\N.
	\end{equation}  
	
\end{theorem}

\begin{remark}
	(a) The distribution of the random index $N$ given in \eqref{eq:sumlog} may be regarded as a mixture with weights
	$\theta$ and $1-\theta$ of the geometric distributions with parameters $1-\theta$ and $\theta$ respectively.
	
	(b) The distributional equation \eqref{eq:reprlog2} may be rewritten as
	\begin{equation}\label{eq:max2}
	     (X+\delta)^2\ =_{\mathcal D} \ \max\{ W_1, W_2,\ldots,W_N\},
	\end{equation}  
	with the same independence structure but where the $W$-variables now have a Weibull distribution. 
    Incidentally, the maximum in~\eqref{eq:reprlog2} resp.~\eqref{eq:max2} is the analogue of the sum in an
	area that is known as tropical arithmetic, see e.g.~\cite{trop}. 
\end{remark}

\subsection{The hyperbolic secant distribution}\label{subsec:hypsec}

The continuous density $f$ of the hyperbolic secant distribution is given by
\begin{equation}\label{eq:densHypSec}
	f(x) \, = \, \frac{1}{\pi} \,  \sech(x) \ \ \biggl(  \, = \, \frac{1}{\pi \cosh(x)} \, 
	            = \, \frac{2}{\pi \, (e^x+e^{-x}) }\biggr), \quad x\in\R.
\end{equation}
It will now be convenient to rescale the noncentrality
parameter via
\begin{equation}\label{eq:delta_to_theta2}
\theta = \theta(\delta) := \frac{2}{(e^{-\delta}+e^\delta)^2}.
\end{equation}   
We note in passing that $\theta(-\delta)=\theta(\delta)$ and that $0<\theta(\delta)<1/2$ 
for all $\delta\not=0$.

\begin{theorem}\label{thm:hypsec}
	Suppose that $X$ has the density function $f$ given in \eqref{eq:densHypSec}. Let $\delta\not=0$
	and let $\theta$ be as in \eqref{eq:delta_to_theta2}. Then the density $f_\delta$ of $|X+\delta|$
	can be written as 
	\begin{equation}\label{eq:reprHypSec}
		f_\delta(x) \; =\; \sum_{k=0}^\infty w_\theta(k)\, g_k(x)\quad \text{for all } x>0,		
	\end{equation}
	where, for all $k\in\N_0$ and $x>0$,
	\begin{align}\label{eq:reprHypSecDef1}
    w_\theta(k) \ &=\ \frac{\binom{2k}{k}}{2^{2k}}\, (1-\theta)^k \, \theta^{1/2},\\ 
    \label{eq:reprHypSecDef2}
    g_k(x) \ &=\ \frac{2^{2k+3/2}}{\pi \, \binom{2k}{k}}\, 
             \frac{e^x+e^{-x}}{e^{2x}+e^{-2x}}\left( \frac{(e^x-e^{-x})^2}{e^{2x}+e^{-2x}}\right )^k.
    \end{align}
\end{theorem}

\begin{remark}\label{remark:HypSec} (a) In \eqref{eq:reprHypSecDef1} and \eqref{eq:reprHypSecDef2}
	the terms have been normalized so that $w_\theta$ is a 
	probability mass function and $g_k$ is a probability density. In fact, the mixing distribution is the 
	negative binomial distribution with parameters $1/2$ and $\theta$, and a random variable $Y$ with density  
	$g_k$ can be obtained from a random variable $Z$ that has the beta distribution with parameters $1/2$ 
	and $k+1/2$ via $Y=\psi(Z)$, where $\psi=\frac{1}{2}\sech^{-1}$ with $\sech^{-1}$ as the inverse of the
        hyperbolic secant function, restricted to 
	the positive half-line. The well-known relationship 
	between the beta and  gamma distributions now implies that we may
	rewrite the mixture representation of the noncentral hyperbolic secant distributions as
	\begin{equation*}
	      |X + \delta|\; =_{\mathcal D}\; \psi\biggl(\frac{Y_0^2}{Y_0^2 +\sum_{j=1}^{2N+1} Y_j^2}\biggr),
	\end{equation*} 
	with $N,Y_0,Y_1,Y_2,\ldots$ independent, $Y_j$ standard normal for all $j\in\N_0$,
	and  the negative binomial distribution with parameters $1/2$ and $\theta$ as the law of $N$; obviously, this implies~\eqref{eq:hypsec0}.	
	
	(b) A similar result holds for the density of $(X+\delta)^2$ instead of $|X+\delta|$. 
	Indeed, writing temporarily $\tilde f_\theta$ for the density of $(X+\theta)^2$, we get
	\begin{equation*}
		\tilde f_\theta(x) \; = \; \frac{1}{2\sqrt{x}}f_\theta(\sqrt{x}) \; =\; \sum_{k=0}^\infty w_\theta(k)\tilde g_k(x)\ 
			           \text{ for all } x>0,
	\end{equation*}
        with $\tilde g_k(x):= g_k(\sqrt{x})/(2\sqrt{x})$ for all  $k\in\N_0$, $x > 0$.

	(c) In the standard literature on probability and statistics the hyperbolic secant distribution is
	not as widely discussed as other distributions that are symmetric about zero, such as the normal, 
	the logistic, the two-sided exponential, or the Cauchy distribution, for example. For various fascinating 
	connections of the hyperbolic secant distribution with the one-dimensional Brownian motion, planar Brownian 
	motion, and the three-dimensional Bessel process the interested reader is referred to \cite{Levy} and \cite{rowi}. 
	As an example, we specially	mention the distributional identity
	\begin{equation}
		\frac{2}{\pi}X \, =_{\mathcal D}\, \frac{Z}{\sup_{0\le t\le 1}\abs{W(t)}},
	\end{equation}
	where $W=(W(t),t\ge 0)$ is a standard Brownian motion and $Z$ is a standard normal random variable
	independent of $W$.  Some statistical applications of the hyperbolic secant distribution
	are given in \cite{Ding}.
\end{remark}

\subsection{Alternative representations in the normal case}\label{subsec:normal}

Here we present a family of mixture representations for the situation where we start with the normal distribution, 
i.e.\ we consider the distribution $\chi^2_1(\delta^2)$ of $(X+\delta)^2$ where $X$ is standard normal. 
In this section we write
$f_\delta$ for the density of $\chi^2_1(\delta^2)$, and $g_k$ for the density associated with the central chi-squared 
distribution with $k$ degrees of freedom, $k\in\N$. 
The classical representation~\eqref{eq:nzchiqua0} may then be written as
\begin{equation}\label{eq:chidens}
	f_\delta \; = \; e^{-\eta} \, \sum_{k=0}^{\infty} \frac{\eta^k}{k!}\, g_{2k+1}
\end{equation}
with $\eta=\delta^2/2$; this equation will appear as a special case of Theorem~\ref{thm:altnormal} below.
 
The new representation families are indexed by two sequences $p=(p_n)_{n\in\N_0}$ and $q=(q_n)_{n\in\N_0}$ of non-negative 
real numbers where we assume that the corresponding power series
\begin{equation*}
u(t) = \sum_{n=1}^\infty p_n t^n,\quad v(t)=\sum_{n=0}^\infty q_n t^n,
\end{equation*}
converge to finite values on some nonempty interval $[0,t_0)$, with $t_0\le\infty$, and where $p_0,t_0,u,v$ satisfy
\begin{equation*}
	p_0=0, \quad v\not\equiv 0, \quad \lim_{t\uparrow t_0} u(t)=\infty.
\end{equation*}
The sequences $p$ and $q$ can be used to obtain 
two families of discrete distributions via `exponential tilting': For each $\theta\in (0,t_0)$,
\begin{equation*}
	p_\theta(n) := \frac{1}{u(\theta)} \, p_n\theta^n,\  q_\theta(n) := \frac{1}{v(\theta)} \, q_n\theta^n,
	     \quad n\in\N_0, 
\end{equation*}
are probability mass functions, and the associated probability generating functions are given by
\begin{equation*}
	t\mapsto \frac{u(\theta t)}{u(\theta)}, \ 	t\mapsto \frac{v(\theta t)}{v(\theta)},\quad 0\le t\le 1.
\end{equation*}
The assumptions on $p$ imply that $u: [0,t_0)\to [0,\infty)$ is a bijection,
which we may use to transform the noncentrality parameter via
\begin{equation}\label{eq:delta2theta3}
	\theta := u^{-1}(|\delta|). 
\end{equation}
In order to define the family of densities for the mixture components and 
the parametric family of mixture distributions we need some more definitions.
We first note that, for all $k\in\N$,  
\begin{equation*}
	u^k(t)= \sum_{n=k}^\infty p_n^{*k}t^n\quad\text{with} \quad 
	    p_n^{*k}=\sum_{\substack{(j_1,\dots,j_k)\in\N^k,\\j_1+\dots +j_k=n}} p_{j_1}\cdots p_{j_k}.
\end{equation*}
We put  $p_0^{*0}:=1$ and $p_j^{*0}:= 0$ for all $j\in\N$, so that this also holds for $k=0$. Further let
\begin{align}\label{eq:defaltnormal2a}
     b_n \ &:=\ \sum_{k=0}^{\lfloor n/2 \rfloor }
         \left (\sum_{\ell=2k}^{n}p_\ell^{*2k}q_{n-\ell}\right )\frac{1}{k!2^k},\quad n\in\N_0,\\
     \label{eq:defaltnormal2}
     B \ &:=\ \{n\in\N_0:\, b_n>0\},\\ 
     \label{eq:defaltnormal3}
     h_n(x)\ &:=\ b_n^{-1}\sum_{k=0}^{\lfloor n/2 \rfloor }\left (\sum_{\ell=2k}^{n}p_\ell^{*2k}q_{n-\ell}\right )
                                          \frac{1}{k!2^k}\,g_{2k+1}(x),\quad x>0,\ n\in B,\\
     \label{eq:defaltnormal4}
     w_{\theta}(n)\ &:=\ \left (v(\theta)\right )^{-1}\exp\left (-u^2(\theta)/2\right )\, b_n \theta^n, \quad n\in \N_0.    
\end{align}
It is clear that the functions $h_n$, $n\in B$, are probability densities on $(0,\infty)$. 
Further, it is part of the assertion of the next theorem that, for each $\theta\in (0,t_0)$, $w_\theta$ is a probability mass function with support $B$.

\begin{theorem}\label{thm:altnormal}
	Let $\theta$, $B$, $h_n$ and $w_\theta$ be as in \eqref{eq:delta2theta3}, \eqref{eq:defaltnormal2},
	\eqref{eq:defaltnormal3} and \eqref{eq:defaltnormal4} respectively. For each $n\in B$ let $Y_n$
	be a random variable with density $h_n$ and let $N$  be a $B$-valued random variable with 
	probability mass function $w_\theta$ that is independent of the $Y$-sequence. Then 
	\begin{equation}\label{eq:repraltnormal}
		(X+\delta)^2\; =_{\mathcal D}\; Y_N	.
	\end{equation}  
\end{theorem}

\begin{remark}\label{rem:altnormal}
	As in the logistic and the hyperbolic secant case, more can be said about the 
	distribution of the random summation index $N$ in \eqref{eq:repraltnormal}: It follows from the computations in the
	proof of Theorem~\ref{thm:altnormal} in Section~\ref{subsec:proofaltnormal} that the generating
	function of $N$ is given by
	\begin{equation*}
	    \frac{v(\theta t)}{v(\theta)}\exp\left (\frac{1}{2}u^2(\theta)\left [\frac{u^2(\theta t)}{u^2(\theta)}-1\right ] \right ),
	    \quad 0\le t\le 1.	
	\end{equation*}
	Based on this explicit formula standard calculations show that 
	\begin{equation}\label{eq:repraltnormalN}
	     N \; =_{\mathcal D} \; K\, +\, \sum_{j=0}^M (L_{j,1}+L_{j,2}).   
	\end{equation} 
	Here the $\N_0$-valued random variables $K,M,L_{j,1},L_{j,2}$, $j\in\N$, are independent, 
	$K$ has the generating function ${v(\theta t)}/{v(\theta)}$,  $0\le t\le 1$, 
	$M$ has the Poisson distribution with parameter
	$\delta^2/2$, and the $L_{j,1},L_{j,2}$, $j\in\N$,
	have the distribution with generating function $u(\theta t)/u(\theta)$, $0\le t\le 1$; see also
	the above remarks on exponential tilting.

\end{remark}

Theorem~\ref{thm:altnormal} is quite general; as a result, the formulas for the mixture components and the 
mixture distribution are somewhat involved. Specializing $u$ and $v$ we obtain several interesting subrepresentations
where some of these ingredients can be made more explicit.

\begin{example}\label{ex:altnormal} 
	(a) The classical case \eqref{eq:nzchiqua0} resp.\ \eqref{eq:chidens} appears if we take $v(t)\equiv 1$, $u(t)=t$ for all $t\in [0,t_0)$ with $t_0:=\infty$,
	or, equivalently, $p_1=1$ and $p_j=0$ for $j\not=1$, $q_0=1$ and $q_j=0$ for $j\not=0$. 
	
	(b)  Let $\alpha\in (-1,\infty)$. With $v(\theta)=(1-\theta)^{-(1+\alpha)}$
	and $u(\theta)=\theta/(1-\theta)$, $\theta\in (0,1)$, some simplifications occur. First,
	we have, for $n\ge 2k$, 
	\begin{equation*}
	\sum_{\ell=2k}^{n}p_\ell^{*2k}q_{n-\ell}\, =\, \binom{n+\alpha}{n-2k},
	\end{equation*}
	and the distributions of $K$ and the $L$-variables in \eqref{eq:repraltnormalN} turn out to 
	be negative binomial with parameters $1+\alpha$ and $1-\theta$, and geometric with
	parameter $1-\theta$, respectively.
	
	(c) The following may be regarded as the $\alpha=-1$ version of the previous part. Again we take 
	$u(\theta)=\theta/(1-\theta)$, but now let $v\equiv 1$. We exploit the identity  
	\begin{equation*}
		\exp\left (\frac{x\theta}{1-\theta}\right )\; = \;
		1+\sum_{n=1}^\infty \theta^n\sum_{k=1}^n\frac{x^k}{k!}\binom{n-1}{k-1},\quad x\in\R,\,\theta\in (-1,+1)
	\end{equation*}
	to obtain   $b_{0}\, =\, 1$, $b_{1}=0$, and, for $n\ge 2$,
	\begin{equation}\label{eq:bsum}
	     b_{n}\ =\ \sum_{k=1}^{\lfloor n/2\rfloor}\binom{n-1}{2k-1}\frac{1}{2^kk!}.
	\end{equation}
	Further, recalling that $g_k$ denotes the density of $\chi_k^2$ for all $k\in \N$, $h_{0}(x)=g_1(x)$ and, 
	for $n\ge 2$ and $x>0$, 
	\begin{equation}\label{eq:hsum}
 	     h_{n}(x)\ =\ b_{n}^{-1}\sum_{k=1}^{\lfloor n/2\rfloor}
 	                              \binom{n-1}{2k-1}\frac{1}{2^kk!} \, g_{2k+1}(x).
	\end{equation}
	In order to obtain a closed-form expression for the sums appearing in \eqref{eq:bsum} and \eqref{eq:hsum}  
	and/or a probabilistic interpretation of the mixture distribution with the density
	$h_{n}$ for $n\ge 2$ it may be helpful to write
	\begin{equation*}
	       b_n \; =\; \sum_{k=1}^{\lfloor n/2\rfloor}L(n,2k)\frac{(2k-1)!!}{n!},\quad n\ge 2,
	\end{equation*}
	and 
	\begin{equation*}
	     h_{n}(x)=b_{n}^{-1}\sum_{k=1}^{\lfloor n/2\rfloor}L(n,2k)\frac{(2k-1)!!}{n!}\,g_{2k+1}(x),\quad x>0,\; n\ge 2,
	\end{equation*}
	where 
	\begin{equation*}
		L(n,\ell)=\frac{n!}{\ell!}\binom{n-1}{\ell-1},\quad n,\ell\in\N,\; \ell\le n,
	\end{equation*}
	are the special (unsigned) Lah numbers: 
	$L(n,\ell)$ is the number of ways a set of $n$ elements can be partitioned into $\ell$
	non-empty linearly ordered subsets. Additionally, $L(0,0)=0$, $L(n,0)=0$ for
	$n>0,$ and $L(n,k)=0$ for $k>n$;  see \cite{Lah} or, e.g., \cite{cer} for recent research on
	Lah numbers.
	
	(d) With $v\equiv 1$ and $u(\theta)=-\log (1-\theta)$, $\theta\in (0,1)$, 
	\begin{equation*}
	   \frac{u(\theta t)}{u(\theta)}\; =\; \frac{\log (1-\theta t)}{\log (1-\theta)},\quad 0\le t\le 1.
	\end{equation*}
	This is the generating function of the logarithmic series distribution with parameter $\theta$, with
	probability mass function $n^{-1}\theta^n/\left (-\log (1-\theta)\right )$, $n\in\N$. 
	Because of
	\begin{equation*}
	\frac{1}{2}\left (\log (1-x)\right )^2\; =\; \sum_{n=2}^\infty \frac{H_{n-1}}{n}x^n,\quad \abs{x}<1,
	\end{equation*}
	where $H_{m}$ the $m$th harmonic number $H_m=\sum_{j=1}^m\frac{1}{j}$, $m\in\N$,  
	see~e.g.\ formula 1.516 in \cite{gr}, we get 
	\begin{equation*}
	    \left (\frac{\log(1-\theta t)}{\log (1-\theta)}\right )^2\; =\; \sum_{n=2}^\infty 2\frac{H_{n-1}}{n}
		       \frac{1}{\left (\log (1-\theta)\right )^2}\theta^n t^n,\quad 0\le t\le 1,
	\end{equation*}
	as the generating function of the distribution of $L_{j,1}+L_{j,2}$, $j\in\N$. 
\end{example}

The family $h_n$, $n\in B$, of probability densities in our alternative representations for 
noncentral chi-squared distributions consists of mixtures of central chi-squared distributions with odd
degrees of freedom. The transition from the classical representation~\eqref{eq:chimix} may thus 
be regarded as a linear change of the mixture basis. We now relate such changes to non-canonical mixture representations
of the Poisson family.

Of course, any distribution with support $\N_0$ can canonically be seen as a mixture of the one-point masses 
concentrated at $k$, $k\in\N_0$. More general mixture representations of the Poisson distributions
consist of two families of distributions on $\N_0$, which we specify by their probability mass functions
$\{c_n:\, n\in\N_0\}$ and $\{m_\eta:\, \eta > 0\}$, requiring 
\begin{equation}\label{eq:mixPo}
      e^{-\eta}\frac{\eta^k}{k!} \; = \; \sum_{n=0}^\infty m_\eta(n)\, c_n(k)\quad \text{for all }k\in\N_0,\, \eta>0. 
\end{equation}
If $B_n$, $n\in\N_0$, are random variables with $P(B_n=k)=c_n(k)$ for all $k\in\N_0$ and $M_\eta$ is a random variable,
independent of the $B$-family, with $P(M_\eta=n)=m_\eta(n)$ for all $n\in\N_0$, then \eqref{eq:mixPo} implies that $N:=B_{M_\eta}$ is Poisson distributed
with parameter $\eta$. Further, starting with the classical representation $(X+\delta)^2=Y_N$,
see~\eqref{eq:repr3} with $\delta^2/2=\eta$, assuming the independence of the random variables $M_\eta,B_n,n\in\N_0,Y_k,k\in \N_0,$
we then obtain $(X+\delta)^2=Y_{B_{M_\eta}}$, which can be rewritten as
\begin{equation*}
	(X+\delta)^2\, =\, Z_N, \quad Z_j:=Y_{B_j}\ \text{ for all } j\in \N_0.
\end{equation*}
This shows that any mixture representation of the Poisson family leads to a mixture representation 
of the noncentral chi-squared distributions. 

Interestingly, there is a converse.

\begin{theorem}\label{thm:mixPo}
	Let $b_n$ and $B$ be as in~\eqref{eq:defaltnormal2a} and~\eqref{eq:defaltnormal2} respectively.
	For $n\in B$ and $k\in\N_0$ put
	\begin{equation*}
        c_n(k)\, := \, \frac{1}{b_nk!2^k}\, \sum_{\ell=2k}^n p_l^{*2k}q_{n-l}\quad\text{if}~n\ge 2k,~\text{and}~c_n(k)=0~\text{otherwise.}
	\end{equation*}
        For $n\in\N_0, n\notin B$, let $c_n$ be some arbitrary probability mass function on $\N_0.$  
        Further, let $\theta$ and 
	$w_\theta(n)$ be as in~\eqref{eq:delta2theta3} and~\eqref{eq:defaltnormal4} respectively.
        Put $\eta=u^2(\theta)/2$ and define the probability mass function $m_\eta$ by
        $m_\eta(n)=w_\theta(n),\,n\in\N_0.$   
        Then the families $\{ c_n:\, n\in \N_0\}$ and $\{m_\eta:\, \eta>0\}$ provide a mixture 
	representation of the Poisson family.
\end{theorem}

\section{Proofs}\label{sec:proofs}

Before we begin with the proofs of the theorems we briefly recall a specific argument for the
classical case as it also underlies our approach in the present paper.

With $f_\delta$ the density of $(X+\delta)^2$, $X$ standard normal, we get
\begin{align*}
	f_\delta(x)\ &=\ \frac{d}{dx}\left(P(X\le \sqrt{x}-\delta) - P(X\le -\sqrt{x}-\delta)\right) \\
                      &=\ \frac{1}{2}\frac{1}{\sqrt{2\pi x}} \left(e^{-(\sqrt{x}-\delta)^2/2}
                      	                                          + e^{-(-\sqrt{x}-\delta)^2/2}\right)	 \\
                      &=\ \frac{1}{\sqrt{2\pi x}}  \, e^{-\delta^2/2} \, e^{-x/2} \, 
                                 \frac{1}{2}\left(e^{\delta \sqrt{x}} + e^{-\delta \sqrt{x}}\right) \\
                       &=\ \sum_{n=0}^{\infty} e^{-\delta^2/2}\left(\frac{\delta^2}{2}\right)^n\; 
                               \frac{2^n}{\sqrt{2\pi}(2n)!} \, x^{n-\frac{1}{2}} e^{-x/2}.             
 \end{align*}
Using the duplication formula for the gamma function, 
\begin{equation*}
       n!\, 2^{2n}\, \Gamma\left(n +\frac{1}{2}\right)\; =\; \pi^{1/2}\, (2n)!,\quad n\in\N_0,
\end{equation*}
this is easily identified as a Poisson mixture with parameter $\delta^2/2$ of the densities
\begin{equation}\label{eq:densChi}
	g_{2n+1}(x) \; =\; \frac{1}{2^{n+\frac{1}{2}}\Gamma\left(n +\frac{1}{2}\right)}\, x^{n-\frac{1}{2}}\, e^{-x/2}\quad x>0,
\end{equation}
 of the distributions $\chi^2_{2n+1}$, $n\in\N_0$, and \eqref{eq:nzchiqua0} follows.
 
 In order to indicate how this can be extended to obtain a proof of  the representation mentioned in
 Example~\ref{ex:altnormal}\,(b) we replace the standard series representation of the exponential function
 used above by
 \begin{equation*}
 (1-\theta)^{-1-\alpha} \exp\left(\frac{x\theta}{1-\theta}\right)
                \;=\; \sum_{n=0}^\infty L_n^\alpha(-x) \,\theta^n   
 \end{equation*}
 for all $\theta\in (-1,1)$ and all $x\in\R$, where
 \begin{equation*}
 	L_n^\alpha(x)\; := \; \sum_{k=0}^{n} \binom{n+\alpha}{n-k}\frac{1}{k!}\, (-x)^k
 \end{equation*}
 are the generalized Laguerre polynomials; see e.g.\ \cite{Erd}. With $\theta:=\abs{\delta}/(1+\abs{\delta})$
 we  get
 \begin{equation*}
 	e^{\delta\sqrt{x}}+e^{-\delta\sqrt{x}}\; 
 	                   =\; \exp\left(\frac{\theta\sqrt{x}}{1-\theta}\right)+ \exp\left(-\frac{\theta\sqrt{x}}{1-\theta}\right),
 \end{equation*}
 so that again the odd terms in the expansion disappear, and a brute force approach finally yields the 
 desired representation.
 
\subsection{Proof of Theorem~\ref{thm:logistic}}
For $\delta\not=0 $ the distribution function $F_\delta$ of $|X+\delta|$ is given by
\begin{equation*}
     F_\delta(x)\ =\ F(x-\delta)-F(-x-\delta)\ 
                  =\ \frac{1}{1+e^{-x+\delta}}-\frac{1}{1+e^{x+\delta}},\quad x\ge 0.
\end{equation*}
With $\theta$ as in \eqref{eq:delta_to_theta} this leads to
\begin{align*}
     \frac{1}{1+e^{-x+\delta}}\ &=\ \frac{\frac{1}{1+e^\delta}}{1-\frac{e^\delta}{1+e^\delta}(1-e^{-x})}\\
                                  &=\ \frac{1-\theta}{1-\theta (1-e^{-x})}\\
                                  &=\ \sum_{n=0}^\infty (1-\theta)\theta^n (1-e^{-x})^n\\
                                  &=\ 1-\theta+\sum_{n=1}^\infty (1-\theta)\theta^n (1-e^{-x})^n,
\end{align*}
and similarly,
\begin{equation*}
     \frac{1}{1+e^{x+\delta}}\ =\ 1-\theta-\sum_{n=1}^\infty \theta (1-\theta)^n(1-e^{-x})^n,
\end{equation*}
both for all $x\ge 0$. Combining these we obtain
\begin{equation}\label{eq:logis}
      F_\delta(x)\ =\ \sum_{n=1}^\infty \bigl((1-\theta)\theta^n+\theta (1-\theta)^n\bigr) (1-e^{-x})^n \quad\text{for all }x\ge 0.
\end{equation}
Now let $E_1,E_2,\dots$ be independent random variables, where for each $j\in\N$ the random variable
$E_j$ has the exponential distribution with mean 1.    
Then, for each $n\in\N$
\begin{equation}
	P(\max(E_1,\dots,E_n)\le x)\, = \, (1-e^{-x})^n,\ x\ge 0,
\end{equation}
and~\eqref{eq:logis} translates into
\begin{equation}\label{eq:logis2}
\abs{X+\delta}\ =_{\mathcal D} \ \max(1,\dots,E_N)
\end{equation}
with $N$ as in the theorem and independent of $E_j$, $j\in\N$.
By the R\'enyi-Sukhatme representation, see e.g.\ p.~721 in \cite{ShoWe},
\begin{equation*}
	\max(E_1,\dots,E_n)\; =_{\mathcal D}\; \sum_{j=1}^n\frac{1}{j}E_j\quad\text{for all }n\in\N.
\end{equation*}
Thus, using \eqref{eq:logis2},
\begin{equation*}
	\abs{X+\delta}
	              \ =_{\mathcal D} \ \sum_{j=1}^N\frac{1}{j}E_j.
\end{equation*}
Also, \eqref{eq:logis2} immediately leads to the representation of the distribution of $(X+\delta)^2$ as a maximum
\begin{equation*}
	(X+\delta)^2 \, =_{\mathcal D}\, \max(E_1^2,\dots,E_N^2)
\end{equation*}
of a random number of squared exponentials. The latter are easily seen to have the Weibull distribution 
with distribution function
$1-e^{-\sqrt{x}}$, $x\ge 0$. 

\subsection{Proof of Theorem~\ref{thm:hypsec}}
With $f$ as in \eqref{eq:densHypSec} and $\delta\not= 0$ the density $f_\delta$ of $|X+\delta|$ is given by
\begin{align*}
       f_\delta(x)\ &=\ f(x-\delta)+f(x+\delta)\\
             &=\ \frac{2}{\pi}\frac{(e^x+e^{-x})(e^\delta+e^{-\delta})}{e^{2x}+e^{-2x}+e^{2\delta}+e^{-2\delta}}\\  
             &=\ \frac{2}{\pi}
                                \frac{(e^x+e^{-x})(e^\delta+e^{-\delta})}{(e^{2x}+e^{-2x})(e^{2\delta}+e^{-2\delta})}
                                \frac{2}{\frac{2}{e^{2x}+e^{-2x}}+\frac{2}{e^{2\delta}+e^{-2\delta}}}\\ 
             &=\ \frac{2}{\pi}
                       \frac{(e^x+e^{-x})(e^\delta+e^{-\delta})}{(e^{2x}+e^{-2x})(e^{2\delta}+e^{-2\delta})}
                       \frac{2}{1+\frac{2}{e^{2\delta}+e^{-2\delta}}-\left (1-\frac{2}{e^{2x}+e^{-2x}}\right )}\\
             &=\ \frac{4}{\pi}
                       \frac{(e^x+e^{-x})(e^\delta+e^{-\delta})}{(e^{2x}+e^{-2x})(e^{2\delta}+e^{-2\delta})}
                       \frac{1}{1+\frac{2}{e^{2\delta}+e^{-2\delta}}}
                       \left (1-\frac{1-\frac{2}{e^{2x}+e^{-2x}}}{1+\frac{2}{e^{2\delta}+e^{-2\delta}}}\right )^{-1}\\
             &=\ \frac{4}{\pi}
                       \frac{e^x+e^{-x}}{e^{2x}+e^{-2x}}\frac{1}{e^\delta+e^{-\delta}}
                       \sum_{k=0}^\infty \left( \frac{(e^x-e^{-x})^2}{e^{2x}+e^{-2x}}\right )^k
                       \left( \frac{e^{2\delta}+e^{-2\delta}}{(e^{\delta}+e^{-\delta})^2}\right )^k,\quad x>0.
\end{align*}
With $\theta$ as in \eqref{eq:delta_to_theta2} this can be written as 
\begin{equation}\label{eq:repr_fdelta}
       f_\delta(x)\ =\ \sum_{k=0}^\infty \frac{2^{3/2}}{\pi}\frac{e^x+e^{-x}}{e^{2x}+e^{-2x}}
                                               \left(\frac{(e^x-e^{-x})^2}{e^{2x}+e^{-2x}}\right )^k
                                                                   (1-\theta)^k\theta^{1/2},\,x>0.
\end{equation}
For each $k\in\N_0$ let
\begin{equation}
    a_k\; :=\; \frac{2^{3/2}}{\pi}\int_0^\infty \frac{e^x+e^{-x}}{e^{2x}+e^{-2x}}
                         \left( \frac{(e^x-e^{-x})^2}{e^{2x}+e^{-2x}}\right )^k\,dx.
\end{equation}
Because of $\int f_\delta(x)\, dx=1$ we have $1=\sum_{k=0}^\infty a_k (1-\theta)^k\theta^{1/2}$,
or, equivalently,
\begin{equation*}
	\sum_{k=0}^\infty a_k (1-\theta)^k\ =\ \theta^{-1/2}\ =\ (1-(1-\theta))^{-1/2}\ 
	      =\ \sum_{k=0}^\infty(-1)^k\binom{-1/2}{k}(1-\theta)^k,
\end{equation*}
from which we deduce that 
\begin{equation*}
    a_k\, =\, (-1)^k\binom{-1/2}{k}\, =\, 2^{-2k}\binom{2k}{k}\quad \text{for all }k\in\N_0.
\end{equation*}
The assertion \eqref{eq:reprHypSec} of the theorem now follows from \eqref{eq:repr_fdelta} and the 
definitions~\eqref{eq:reprHypSecDef1} and~\eqref{eq:reprHypSecDef2}.

The proof of the distributional fact noted in Remark \ref{remark:HypSec} (a) that $\frac{1}{2}\sech^{-1}(Y)$ has the density $g_k$
given in~\eqref{eq:reprHypSecDef2}, if $Y$ has the beta distribution with parameters $1/2$ and $k+1/2,$
is easily carried out by straightforward calculation. 

\subsection{Proof of Theorem~\ref{thm:altnormal}}
\label{subsec:proofaltnormal}
Again, we begin with a suitable calculation, 
\begin{align*}
     v(\theta)\exp\left (xu(\theta)\right )
     \ &=\ \sum_{n=0}^\infty q_n\theta^n\ +\ \left (\sum_{\ell=0}^\infty q_\ell \theta^\ell \right)\;
         \left(\sum_{k=1}^{\infty}\frac{1}{k!}x^ku(\theta)^k\right ) \\
       &=\ \sum_{n=0}^\infty q_n\theta^n\ +\ \left (\sum_{\ell=0}^\infty q_\ell \theta^\ell \right)\;
         \left(\sum_{k=1}^{\infty}\frac{1}{k!}x^k\sum_{\ell=k}^\infty p_\ell^{*k}\theta^\ell\right)\\
       &=\ \sum_{n=0}^\infty q_n\theta^n\ +\ \left (\sum_{\ell=0}^\infty q_\ell \theta^\ell \right)\;
         \left(\sum_{\ell=1}^{\infty}\left [\sum_{k=1}^\ell\frac{1}{k!}p_\ell^{*k}x^k\right ]\theta^\ell\right)\\
       &=\ \sum_{n=0}^\infty \left(q_n\theta^n+ \sum_{\ell=1}^n \left(\sum_{k=1}^\ell 
                                   \frac{x^k}{k!}p_\ell^{*k}q_{n-\ell}\right )\theta^n\right ).  
\end{align*}
With
\begin{equation*}
	L_n(x)\; :=\; \sum_{\ell=0}^{n}\left (\sum_{k=0}^\ell\frac{x^k}{k!}p_\ell^{*k}\right )q_{n-\ell}
	\; =\; \sum_{k=0}^n \left (\sum_{\ell=k}^n p_\ell^{*k} q_{n-\ell}\right )\frac{x^k}{k!}.
\end{equation*}
this may be rewritten as
\begin{equation*}
	v(\theta)\exp\left (xu(\theta)\right )\, =\, \sum_{n=0}^\infty L_n(x)\, \theta^n.
\end{equation*}
If $q_0\neq 0$ and $p_1\neq 0$ the polynomials $L_n(x)$, $n\in\N_0$, are of the generalized Appell type;
more specifically, the $n!L_n(x)$ arise as Sheffer polynomials, see~\cite{roman}. Further, with
\begin{equation*}
   M_n(x)\; :=\; \sum_{k=0}^{\lfloor n/2\rfloor } \left (\sum_{\ell=2k}^{n} p_\ell^{*2k} q_{n-\ell}\right )\frac{x^{k}}{(2k)!}	
\end{equation*}
we get 
\begin{equation*}
	M_n(\sqrt{x})=\frac{1}{2}\left[L_n(\sqrt{x})+L_n(-\sqrt{x})\right ]\quad \text{for }x>0
\end{equation*}
and therefore   
\begin{equation*}
    \frac{1}{2}\left (\exp(\sqrt{x}\delta)+\exp(-\sqrt{x}\delta)\right) \ 
                   =\ \frac{1}{v(\theta)}\sum_{n=0}^{\infty} M_n(x)\, \theta^n.
\end{equation*}
Using the definitions in \eqref{eq:defaltnormal2a} to \eqref{eq:defaltnormal4} we may now write the density
$f_\delta$ of $\chi^2_1(\delta^2)$ as
\begin{align}\label{eq:allgdar}
    f_\delta(x)=\sum_{n\in B} h_n(x)w_\theta(n),\quad x>0.
\end{align}
It is clear that $h_n$ is a probability density for all $n\in  B$. Finally, to see that $w_\theta$ is a probability 
mass function, we first note that $w_\theta(n)\ge 0$ and then use \eqref{eq:allgdar} to obtain 
\begin{align*}
   1\, =\, \int_0^\infty f_\delta(x)\,dx\, =\, \sum_{n\in B}w_\theta(n).
\end{align*}

\subsection{Proof of Theorem~\ref{thm:mixPo}}

Taken together, and with $g_k$ the density of $\chi^2_k$, the classical 
representation~\eqref{eq:chidens}
and the alternative representation given in Theorem~\ref{thm:altnormal} imply that
\begin{equation*}
	\sum_{k=0}^{\infty} e^{-\eta}\frac{\eta^k}{k!}\, g_{2k+1} \ = \ f_\delta
	 \ = \ \sum_{n=0}^{\infty}\sum_{k=0}^{\infty}m_\eta(n)c_n(k)\, g_{2k+1}.
\end{equation*}
As all the individual terms are nonnegative we may rearrange the right hand side. 
Multiplying both sides by $\sqrt{2x}e^{x/2}$ we then 
obtain
\begin{equation*}
	\sum_{k=0}^{\infty} e^{-\eta}\frac{\eta^k}{k!}\, \frac{1}{2^k\Gamma\left(k+\frac{1}{2}\right)}\, x^k
	  \ =\ \sum_{k=0}^{\infty}\biggl( \sum_{n=0}^\infty m_\eta(n) c_n(k) \biggr)
	                   \frac{1}{2^k\Gamma\left(k+\frac{1}{2}\right)}\, x^k,
\end{equation*}
and
\begin{equation*}
	e^{-\eta}\frac{\eta^k}{k!}\ = \ \sum_{n=0}^\infty m_\eta(n) c_n(k)
\end{equation*}
follows by comparison of coefficients.

\section{Comments}\label{sec:comments}

\subsection{The connection to statistics}\label{subsec:stat}

So far the only noncentral chi-squared distributions considered were those with one degree of freedom.
In statistics, the distribution of $Z:= \sum_{i=1}^k (X_i+\mu_i)^2$ is important, where $X_1,\ldots,X_k$ are independent
standard normals and $\mu_1,\ldots,\mu_k$ are arbitrary real numbers.
If $\delta^2:=\sum_{i=1}^k \mu_i^2>0$ this is the noncentral
chi-squared distribution $\chi^2_k(\delta^2)$ with $k$ degrees of freedom and noncentrality parameter $\delta^2$,
and the following generalization 
of~\eqref{eq:chimix} holds:
\begin{equation*}\label{eq:chigen_dof}
\chi^2_k(\delta^2) \;=\; e^{-\delta^2/2} \sum_{n=0}^\infty \frac{\delta^{2n}}{2^n n!}\, \chi^2_{2n+k}
\quad\text{for all } \delta\not=0, \, k\in\N.
\end{equation*}
In contrast to this, the known mixture representations for noncentral $F$-distributions,  
see e.g.\ \cite{stuartord},
do not arise as the distributions of a function of $|X+\delta|$, $\delta\not=0$, for some random variable $X$ satisfying $X=_{\mathcal D} -X$.
Noncentral $t$-distributions are not even concentrated on the nonnegative half-line, 
but see~\cite{stuartord} again and Section~\ref{subsec:geo} below.

Special quadratic forms of multivariate normal random vectors have noncentral chi-squared distributions;
see, e.g., \cite{rao} for a thorough discussion. Additionally, the power function of various
important statistical tests can be expressed in terms of noncentral chi-squared distributions.
To give a simple specific example,
let $X$ be a $d$-variate normal random vector with unknown mean vector $a\in\R^d$ and the $d\times d$
identity matrix as covariance matrix. Consider testing the hypothesis $H:a=0$ against the general alternative
$K:a\neq 0.$ The testing problem is invariant under the group of orthogonal transformations on $\R^d.$
For given $\alpha\in (0,1),$ there exists a uniformly most powerful invariant test at level $\alpha.$ With the
euclidean norm $\norm{X}^2$ of $X$ there is a simple maximal invariant test statistic.
Denoting by $\chi_{d;1-\alpha}^2$ the $1-\alpha$ quantile of $\chi_d^2,$ the test rejects the hypothesis iff
$\norm{X}^2>\chi_{d;1-\alpha}^2.$ If $X$ has the mean $a\neq 0$ the distribution of the
test statistic $\norm{X}^2$ is the noncentral $\chi_d^2$ distribution with noncentrality
parameter $\norm{a}^2.$ Thus, denoting by $G_{d,\delta^2}$ the distribution function of the $\chi_d^2(\delta^2)$
distribution, the power function of the test is given by $1-G_{d,\norm{a}^2}(\chi_{d;1-\alpha}^2),a\in \R^d.$
See, e.g., \cite{liese} for other testing problems where the power function can be
expressed in terms of noncentral chi-squared distributions.

\subsection{The connection to stochastic processes}\label{subsec:stochproc}
Let $X=(X_t)_{t\ge 0}$ be an irreducible Markov chain with finite state space $S$ and
symmetric jump rates. For $i\in S$ and $t\ge 0$ let
\begin{equation*}
	L^0_i(t) \, := \, \int_0^t 1\{X_u=i\}\, du
\end{equation*} 
be the amount of time spent by $X$ in state $i$ up to time $t$. For a fixed state $i_0$ 
and a given $\eta>0$ let
\begin{equation*}
	\tau \,:=\, \inf\{t\ge 0: L_{i_0}(t) \ge \eta\},
\end{equation*}
and let $L=(L_i)_{i\in S'}$ with
$S':=S\setminus\{i_0\}$
be given by $L_i:=L_i^0(\tau)$.
Suppose that $X$ starts at $i_0$.
The symmetry assumption implies that the matrix
\begin{equation*}
	G=(g_{ij})_{i,j\in E'}, \quad g_{ij}:= E L_iL_j,
\end{equation*}
is positive semidefinite and symmetric; in particular, it can serve as the covariance matrix
of a central normal random vector $Y=(Y_i)_{i\in S'}$. The generalized second Ray-Knight 
theorem says that, for $X$ and $Y$ independent,
\begin{equation*}\label{eq:RK}
     \frac{1}{2}Y^2 + L \; =_{\mathcal D}\; \frac{1}{2}\bigl(Y+\sqrt{2\eta})^2.
\end{equation*}
The standard reference for this circle of ideas is the book by~\cite{iso}, but see also the
very recent approach in~\cite{isoneu}.

The simplest nontrivial example has $S=\{0,1\}$ and jump rates 1. We take $i_0=0$
and glue the time intervals in 0 together, which results in a Poisson process with rate 1.
This implies that the number $N$ of visits to $1$ up to time $\tau$ has a Poisson distribution with 
parameter $\eta$. The $L$-vector consists of one component only, it is the sum of
$N$ independent random variables $E_1,E_2,\ldots$ that are all exponentially distributed with mean 1.
The $G$-matrix consists of the single number 1, which means that $Y$ is a standard normal
random variable. Taken together we see that~\eqref{eq:RK} leads to
\begin{equation*}
 \bigl(Y+\sqrt{2\eta})^2\; =_{\mathcal D}\; Y^2 \, +\, 2\,\sum_{j=1}^N E_j,
\end{equation*} 
which is \eqref{eq:nzchiqua0}.

The classical second Ray-Knight theorem refers to Brownian motion where a discrete time approximation
leads to the simple symmetric random walk. In the approach given by~\cite{moper} the proof is reduced
to the two-states situation and~\eqref{eq:nzchiqua0}. The authors of the book ask for a probabilistic argument 
(as the authors of the present paper did more than 20 years ago, coming from the statistical side).
In fact, in the much more general circle of ideas known as the Dynkin isomorphism, the lack of a
probabilistic explanation of the basic distributional equalities is often mentioned; see e.g.~the preface
to~\cite{iso}.

\subsection{A geometric outlook}\label{subsec:geo}
The mixture representation $\mu=\sum_{k=0}^\infty a_k\nu_k$
of  a probability distribution $\mu$ in terms of a mixing 
distribution $a_k$ and a mixing base $\nu_k$ of probability distributions, $k\in\N_0$,  may 
equivalently be seen 
as a representation of $\mu$ as a convex combination of
the $\nu_k$'s. In particular,~\eqref{eq:nzchiqua0} leads to
\begin{equation}\label{eq:clconv1}
	\{\chi_1^2(\delta):\,\delta\not=0\}\; \subset\; \overline{\text{conv}\{\chi^2_{2k+1}:\, k\in \N_0\}},
\end{equation}
where the notation on the right refers to the closure of the convex hull, and closure in turn refers to 
total variation norm (or $L^1$-convergence if there are densities).  Similarly,~\eqref{eq:chigen_dof} implies
\begin{equation}\label{eq:clconv2}
\{\chi_k^2(\delta):\,\delta\not=0,\, k\in\N\}\; \subset\; \overline{\text{conv}\{\chi^2_k:\, k\in \N\}}.
\end{equation}
Of course, the mixture results go beyond this geometric interpretation as they also provide the 
mixing coefficients $a_k(\delta)$, $k\in\N_0$, for each $\delta\not=0$.

Moving from convex to general linear combinations $\mu=\sum_{k=0}^\infty a_k\nu_k$ we obtain
\begin{equation*}
	\mu\, \in \, \overline{\text{Lin}\{\nu_k:\, k\in \N_0\}},
\end{equation*} 
where the notation now refers to the closed linear span, and the topology may be the $L^1$-
or $L^2$-distance in the density case (we temporarily do not distinguish between measures and
their densities in this case). For example, the generalized Laguerre polynomials with
$\alpha=1/2$ can easily be transformed into an orthonormal basis of the Hilbert space
$L^2:=L^2(\R,\BB_+,\chi_1^2)$
of all measurable functions $f:(0,\infty)\to\R$ satisfying $\int f^2\,d\chi_1^2<\infty$. In view of 
\begin{equation*}
	\frac{d\chi_{2k+1}^2}{d\chi_1^2}(x) \,=\, c_k x^k\quad \text{with } 
	           \ c_k:= 2^{-k}\, \Gamma(1/2)/\Gamma(k+1/2) = 1/(2k-1)!!
\end{equation*}
the central chi-squared densities with degrees $2k+1$, $k\in\N$, are finite linear combinations of these basis 
elements, and vice versa. This implies 
\begin{equation*}
	\overline{\text{Lin}\{\chi_{2k+1}^2:\, k\in \N\}} \, =\,  L^2
\end{equation*}
so that, for example, the
uniform distributions on $[0,\theta]$, $\theta>0$, can be written as infinite linear combinations
of the central chi-squared distributions with odd degrees of freedom, but it is easy to see that they are not contained in their
closed convex hull. Also, the representation for noncentral $t$-distributions given in~\cite{stuartord} is of the general linear
and not of the convex type.

In both cases the representations lead to approximations by finite sums. For the statistician
and before computers became widely available this aspect was of particular importance,
and it generated a sizable literature. For
example,~\cite{Tiku} used Laguerre polynomials in the context of approximating noncentral
chi-squared distributions, but in contrast to our approach the coefficients can be negative.
An advantage of a mixture representation over a general linear representation 
is that the corresponding approximations for the
distribution function are monotone increasing; also, tail estimates for the mixture distribution
immediately lead to error bounds. It is therefore of considerable interest to obtain similar mixture 
representations for other symmetric distributions, such as the Cauchy distribution and the two-sided
exponential distribution. Finally, the convexity aspect raises some interesting questions of its own. 
For example, we may pass to the closed convex hull on the left hand side of both~\eqref{eq:clconv1} 
and~\eqref{eq:clconv2} without destroying the inclusion relation. 
Does this operation lead to both sets being equal? (It is easy to see that this is indeed the case 
for~\eqref{eq:clconv2}).
Is there a general and usable description of the closed convex hull of general 
parametric families?  We postpone these questions to future research work.

\bigskip
{\bf Acknowledgement.} The authors would like to thank an anonymous reviewer for helpful suggestions.

\end{document}